\magnification 1200
\advance\hsize by 2,3 truecm
\def\makefootline{\baselineskip=52pt\line{\the\footline}}
\vsize= 23 true cm
\hsize= 16 true cm
\overfullrule=0mm

\headline={\hfill\tenrm\folio\hfil}
\footline={\hfill}\pageno=1
\def\vvert#1{\Vert#1\Vert}

\newcount\coefftaille \newdimen\taille
\newdimen\htstrut \newdimen\wdstrut
\newdimen\ts \newdimen\tss

\catcode`\@=11
\def\Eqalign#1{\null\,\vcenter{\openup\jot\m@th\ialign{
\strut\hfil$\displaystyle{##}$&$\displaystyle{{}##}$\hfil
&&\quad\strut\hfil$\displaystyle{##}$&$\displaystyle{{}##}$
\hfil\crcr#1\crcr}}\,}  \catcode`@12

\def\fspeciale{\textfont0=\tenrmp%
\scriptfont0=\sevenrmp%
\scriptscriptfont0=\fivermp%
\textfont1=\tenip%
\scriptfont1=\sevenip%
\scriptscriptfont1=\fiveip%
\textfont2=\tensyp%
\scriptfont2=\sevensyp%
\scriptscriptfont2=\fivesyp%
\textfont3=\tenexp%
\scriptfont3=\tenexp%
\scriptscriptfont3=\tenexp%
\textfont\itfam=\tenitp%
\textfont\bffam=\tenbfp%
\textfont\slfam=\tenbfp%
\def\it{\fam\itfam\tenitp}%
\def\bf{\fam\bffam\tenbfp}%
\def\rm{\fam0\tenrmp}%
\def\sl{\fam\slfam\tenslp}%
\normalbaselineskip=12pt%
\multiply \normalbaselineskip by \coefftaille%
\divide \normalbaselineskip by 1000%
\normalbaselines%
\abovedisplayskip=10pt plus 2pt minus 7pt%
\multiply \abovedisplayskip by \coefftaille%
\divide \abovedisplayskip by 1000%
\belowdisplayskip=7pt plus 3pt minus 4pt%
\multiply \belowdisplayskip by \coefftaille%
\divide \belowdisplayskip by 1000%
\setbox\strutbox=\hbox{\vrule height\htstrut depth\wdstrut width 0pt}%
\rm}

\def\vmid#1{\mid\!#1\!\mid}

\def\fle{\rightarrow}

\null\vskip-1cm

\font\sc=cmcsc10

\newdimen\emm 
\def\pmb#1{\emm=0.03em\leavevmode\setbox0=\hbox{#1}
\kern0.901\emm\raise0.434\emm\copy0\kern-\wd0
\kern-0.678\emm\raise0.975\emm\copy0\kern-\wd0
\kern-0.846\emm\raise0.782\emm\copy0\kern-\wd0
\kern-0.377\emm\raise-0.000\emm\copy0\kern-\wd0
\kern0.377\emm\raise-0.782\emm\copy0\kern-\wd0
\kern0.846\emm\raise-0.975\emm\copy0\kern-\wd0
\kern0.678\emm\raise-0.434\emm\copy0\kern-\wd0
\kern\wd0\kern-0.901\emm}

\font\tendb=msbm10
\font\sevendb=msbm7

\newfam\dbfam
\textfont\dbfam=\tendb\scriptfont\dbfam=\sevendb\scriptscriptfont\dbfam=\sevendb
\def\db{\fam\dbfam\tendb}

\def\C{{\db C }}

\def\N{{\db N }}

\def\R{{\db R }}

\def\Z{{\db Z }}


\def\mod{\mathop{\rm mod}\nolimits}

\font\gothique=eufm10
\def\got#1{{\gothique #1}}

\newdimen\margeg \margeg=0pt
\def\bb#1&#2&#3&#4&#5&{\par{\parindent=0pt
    \advance\margeg by 1.1truecm\leftskip=\margeg
    {\everypar{\leftskip=\margeg}\smallbreak\noindent
    \hbox to 0pt{\hss\bf [#1]~~}{\bf #2 - }#3~; {\it #4.}\par\medskip
    #5 }
\medskip}}

\newdimen\margeg \margeg=0pt
\def\bbaa#1&#2&#3&#4&#5&{\par{\parindent=0pt
    \advance\margeg by 1.1truecm\leftskip=\margeg
    {\everypar{\leftskip=\margeg}\smallbreak\noindent
    \hbox to 0pt{\hss [#1]~~}{\pmb{\sc #2} - }#3~; {\it #4.}\par\medskip
    #5 }
\medskip}}

\newdimen\margeg \margeg=0pt
\def\bba#1&#2&#3&#4&#5&{\par{\parindent=0pt
    \advance\margeg by 1.1truecm\leftskip=\margeg
    {\everypar{\leftskip=\margeg}\smallbreak\noindent
    \hbox to 0pt{\hss [#1]~~}{{\sc #2} - }#3~; {\it #4.}\par\medskip
    #5 }
\medskip}}

\def\messages#1{\immediate\write16{#1}}

\def\findem{\vrule height0pt width4pt depth4pt}

\long\def\demA#1{{\parindent=0pt\messages{debut de preuve}\smallbreak
     \advance\margeg by 2truecm \leftskip=\margeg  plus 0pt
     {\everypar{\leftskip =\margeg  plus 0pt}
              \everydisplay{\displaywidth=\hsize
              \advance\displaywidth  by -1truecm
              \displayindent= 1truecm}
     {\bf Proof } -- \enspace #1
      \hfill\findem}\bigbreak}\messages{fin de preuve}}

\def\resp{\mathop{\rm resp}\nolimits}
\def\resp.{\mathop{\rm resp.}\nolimits}
\font\zpeur=cmtex10
\def\euro{{\zpeur C}\hskip-7pt$=$}

\null\vskip1cm

\parindent=0cm

{\bf HOLOMORPHIC CLIFFORDIAN PRODUCT}
\vskip1cm
{\obeylines Guy Laville
UMR 6139, Laboratoire Nicolas ORESME
D\'epartement de Math\'ematiques
Universit\'e de Caen
14032 CAEN Cedex France
{\it glaville@math.unicaen.fr}
}

\vskip2cm
{\bf Abstract.} \ Let $\R_{0,n}$ be the Clifford algebra of the
antieuclidean vector space of dimension $n$. The aim is to built a
function theory analogous to the one in the \ $\C$ \ case.  In the
latter case, the product of two holomorphic functions is holomorphic,
this fact is, of course, of paramount importance. Then it is necessary
to define a product for functions in the Clifford context.
But, non-commutativity is inconciliable with product of 
functions. Here we introduce a product which is commutative and we
 compute some examples explicitely.



\vskip1cm
{\sc 1. -- Introduction}
\bigskip\medskip 
In one complex variable, it is possible to define a product of two
holomorphic functions   $f$   and   $g$   by  $(fg) (z) =
f(z)g(z)$ because this last expression is holomorphic. Here we make
use of commutativity and of Cauchy-Riemann equations which are first
order partial dif\-ferential equations.
But in fact,
there is much more than that. Holomorphy is equivalent of analyticity~:
taking \ $f(z) = \Sigma a_pz^p$ and $g(z) = \Sigma a_q z^q$ \ then
$$(fg) (z) = \sum_n \ \Big( \sum_{p+q=n} \ a_p \ b_q\Big) z^n.$$

\bigskip
We can do the product either in the space of the values or in the
space of the variable and para\-meters. For higher dimensional spaces,
in Clifford analysis, the above two possibilities give two different
results. The first product is useless because if $f(x)$ and $g(x)$
are monogenic [1], [3], or regular [6], or holomorphic
Cliffordian [9], $f(x)g(x)$ \ is not.  In [1],
 F. {\sc Bracks}, R. {\sc Delanghe},  F. {\sc Sommen}
defined the Cauchy Kovalewski product, but it is no so easy
to work with it [8]. The existence of a product is one of the principal
questions in Clifford analysis, see [11] and [13].  In [7] D. {\sc
Hestenes} and G. {\sc Sobczyk} defined the inner product.  In [10] H.
{\sc Malonek} worked with his permutational product. It is related to
Fueter's ideas [5].

\bigskip
The anticommutator \ $\{ a,b\} = 1/2 (ab+ba)$ \ is well known, but
when we have three elements, we get \ $\{ a, \{ b,c\}\}$ \ or \ $\{ \{
a,b\}, c\}$ \ or \ $\{ \{ a,c\}, b\}$.  In several papers \ [12],
[14], \ {\sc F. Sommen} uses the basic fact that the anticommutator
of two vectors is a scalar and hence commutes with all elements. By
the same token here a basic fact is that the anticommutator of two
paravectors is a paravector.

\bigskip
In quantum mechanics other products are defined~: chronological
product, normal order in product.

\bigskip
Notations.
\bigskip
Let $\R_{0,n}$ the Clifford algebra of the real vector space $V$ of
dimension $n$, provided with a quadratic form of negative signature.
Denote by $S$ the set of scalars in $\R_{0,n}$ which can be
identified to $\R$. An element of the vector space $S\oplus V$ is
called a paravector. Let $\{ e_i\}, i=1,\ldots,n$ be an orthonormal
basis of $V$ and let $e_0 = 1$. We have $e_ie_j+e_je_i = -
2\delta_{ij}$ for $1\leq i,j\leq n$. On $S\oplus V$ we have two
quadratic structures~: one with signature 
$+ - \cdots -$, \  the other with signature $++ \cdots +$.  In this
latter case the scalar product is denoted by $(a\mid b)$.  To do
analysis, we take a norm on $S\oplus V$ such that \
$\vvert{ab}~\leq~\vvert{a}~\vvert{b}$.

For any paravector $u$,  we split up the real part $u_0$ and the
vectorial part $\vec u$~:
$$u = u_0 + \vec u.$$

\vskip1cm
{\sc 2. -- Algebraic structure on the paravector space}
\bigskip
{\bf 2.1  Symmetric product}
\bigskip
{\sc Theorem} and {\sc definition 1}.- \ For $\ell\in\N \setminus \{
0\}$  define the multilinear symmetric function 
$$\Eqalign{
&E : ~(S \oplus V)^\ell &&\longrightarrow~ S \oplus V\cr
&\qquad (u_1,\ldots,u_\ell) &&\longrightarrow~ {1\over\ell!} \
\sum_{\sigma\in \hbox{\got S}_\ell} u_{\sigma (1)}
\ldots u_{\sigma (\ell)}\cr}$$
where $\hbox{\got S}_\ell$ is the set of all permutations
of $\{ 1,\ldots,\ell\}$.

\bigskip
Proof.- \  It is obvious that this function is
multilinear and symmetric. To prove that the values are in
$S\oplus V$, we need a lemma, but before stating it, it is useful to
introduce an algorithmic symbol~:
\vskip-0,8cm
$$~\hbox{\euro}~ \prod_{i=1}^\ell u_i ~\hbox{\euro}~:= {1\over \ell
!} \ \sum_{\sigma\in\hbox{\got S}_\ell} u_{\sigma (1)}
\ldots u_{\sigma (\ell)}.\leqno (1)$$

It is easier to work with this than with \ $E
(u_1,\ldots,u_\ell)$.

\bigskip\bigskip
Lemma 1. \vskip-0,8cm
$$\leqalignno{
&~\hbox{\euro}~ \prod_{i=1}^\ell u_i~\hbox{\euro}~=  &(2)\cr 
&{1\over \ell} \
\sum_{i=1}^\ell u_i ~\hbox{\euro}~ \prod_{\scriptstyle
j=1\atop\scriptstyle j\not= i}^\ell u_j ~\hbox{\euro}~=&\cr 
&{1\over \ell} \ \sum_{i=1}^\ell \ \hbox{\euro} \  \prod_{\scriptstyle 
j=1\atop\scriptstyle j\not= i}^\ell u_j \ \hbox{\euro}  \  u_i = &\cr
&{1\over 2\ell}  \Big( \sum_{i=1}^\ell u_i  \ ~\hbox{\euro}~
\prod_{\scriptstyle  j=1\atop\scriptstyle j\not= i}^\ell 
u_j ~\hbox{\euro}~ + \sum_{i=1}^\ell \  ~\hbox{\euro}~
\prod_{\scriptstyle  j=1\atop\scriptstyle
j\not= i}^\ell u_j ~\hbox{\euro}~ \ u_i\Big).&\cr}$$

\bigskip
Proof of (2).- \ The first and second formulas are
factorisations of the symmetric product. The third one is a
mean of these two.

\bigskip
Now, to prove that the values of the function $E$ is in
$S\oplus V$, we use induction on $\ell$.
\bigskip
For \ $\ell = 1$ \ the result is trivial, for \ 
$\ell = 2$, we have 
$$~\hbox{\euro}~ ab~\hbox{\euro}~= \displaystyle{1\over 2} (ab+ba)$$
is in $S\oplus V$. The last formula (2) allows us to finish
the recurrence.

\bigskip\bigskip
Proposition 1.- \ {\it For  $i=1,\ldots,\ell$ and $u_i\in S
\oplus V$}
$$\vvert{~\hbox{\euro}~ \prod_{i=1}^\ell~Êu_i~\hbox{\euro}~}~ \leq~
\vvert{\prod_{i=1}^\ell u_i}.\leqno (3)$$
This follows from the definition.

\bigskip\bigskip
Extension of the symbol.
\medskip
Let \ $\varphi$ \ be a linear function~:
$$\Eqalign{
&(S \oplus V)^k  &&\longrightarrow~ (S \oplus V)^\ell\cr
&(u_1,\ldots,u_k) &&\longrightarrow~ \big( \varphi_1
(u_1,\ldots,u_k),\ldots,\varphi_\ell
(u_1,\ldots,u_k)\big)\cr}$$

then, we define
\vskip-0,8cm
$$~\hbox{\euro}~ \prod_{i=1}^\ell \varphi_i
(u_1,\ldots,u_p)~\hbox{\euro}~:= E\circ\varphi \
(u_1,\ldots,u_p).\leqno (4)$$

\bigskip
Remark 1.
\medskip
It is always possible to restrict the symmetrization to only
$n = \dim V$ factors because if we have $\ell$ paravectors
$u_1,\ldots, u_\ell$ we can take $\vec
u_{i_{1}},\ldots,\vec u_{i_{p}}$, ~$p$ linearly independent
vectors, all paravectors $u_j$ are linear combinations
of $1$ and the $\vec u_{i_{k}}$ and the symmetrization is
on $\vec u_{i_{1}},\ldots, \vec u_{i_{p}}$.

\bigskip
Remark 2.
\medskip
$~\hbox{\euro}~ \displaystyle\prod_{i=1}^k A_i ~\hbox{\euro}~$ \ is well
defined for all \ $A_i \in\R_{0,n}$ \ because $A_i$ are sums
and products of paravectors and we have linearity.

\bigskip
Remark 3.
\medskip
$~\hbox{\euro}~ x~\hbox{\euro}~ yz~\hbox{\euro}~\hbox{\euro}~$ 
makes sense but it is clumsy and it is a pitfall, so we shall avoid using it. In general it is not
equal to \
$~\hbox{\euro}~ xyz~\hbox{\euro}~$.

\bigskip\bigskip
{\bf 2.2  The symmetric algebra of $V$}
\bigskip\medskip
For \ $n=1$, \ $\R_{0,1} = \C$ \ we have a special phenomena. Take \
$\R [X]$ \ the algebra of polynomials in one indeterminate, then \ $\R
[X]/(X^2+1)$ \ is \ $\C$.  But \ $\R [X_1,\ldots,X_n]$, \ the algebra
of polynomials in \ $n$ \ indeterminates is not directly connected
with \ $\R_{0,n}$. This algebra of polynomials is clearly built to do
products.

Inside the \ $\hbox{\euro}$  \ we compute in \ $\R [e_1,\ldots,e_n]$ \
which may be identified with the symmetric algebra (algebra of
symmetric tensors) of the vector space \ $V$.

\bigskip
{\bf 2.3  Examples}
\bigskip\medskip
In the following formulas $a,b,c$ are in $S\oplus V$.
$$\eqalign{
&~\hbox{\euro}~ a~\hbox{\euro}~=Ê a\cr
&~\hbox{\euro}~ ab~\hbox{\euro}~ = \displaystyle{1\over 2} (ab+ba)\cr
&~\hbox{\euro}~ a^2b~\hbox{\euro}~= \displaystyle{1\over 3} (a^2b +
aba + ba^2).\cr}$$

It is important to notice that this is not \
$\displaystyle{1\over 2} ~(a^2b + ba^2)$

$$\eqalign{~\hbox{\euro}~ e_1^2 a~\hbox{\euro}~&= -
\displaystyle{2\over 3} a +{1\over 3}~ e_1 ~a e_1\cr
\pmatrix{p+q\cr q\cr} ~\hbox{\euro}~ a^p b^q ~\hbox{\euro}~&=
{\displaystyle{d^q\over dt^q}\Big| }_{q=0} (a+tb)^{p+q}.\cr}$$

\bigskip
 Here are explicit formulas for \ $~\hbox{\euro}~
\displaystyle\prod_{h=1}^H e_{i_{h}}~\hbox{\euro}~$~:

If \  $H \equiv 0 ~(\mod~4)$ \  and if all indices are equal then it is equal to $1$, otherwise it is
$0$.

\medskip
If \ $H \equiv 1 ~(\mod~4)$ \  and if all indices are equal then  it is equal to
 $e_{i_{1}}$, 
otherwise if all indices are equal but one, say \  $i_1$,   it is
$\displaystyle{1\over H} e_{i_{1}}$ else it is ~$0$.

\medskip
If \ $H \equiv 2 \ (\mod~4)$ \  and if all indices are equal then it is equal to \  $-1$ \ 
otherwise \  $0$.

\medskip
If \ $H \equiv 3 \ (\mod~4)$ \  and if all indices are equal then  it is equal to \ 
$- e_{i_{1}}$ \ 
otherwise if all indices are equal but one, say \  $i_1$, it is \  $-
\displaystyle{1\over H} e_{i_{1}}$ \   else it is \ $0$.

\bigskip\bigskip
Proof of these values~: 
\smallskip
Take \ $v = (t_1~e_{i_{1}} + \ldots + t_H~e_{i_{H}})$ \ where \
$t_1,\ldots,t_H$ \ are scalars.

Beside the coefficient, the value of the product is the homogeneous
term corresponding to \ $t_1t_2 \ldots t_H$ \ in \ $v^H$.

\medskip
First case~: all indices are equal, say \ $e_{i_{1}}$
$$\Eqalign{
&v   &&\hskip-0,3cm = (t_1 +\ldots+ t_H) e_{i_{1}}\cr
&v^H &&\hskip-0,3cm = (t_1 + \ldots + t_H)^H \ e_{i_{1}}^H\cr}$$
and \ $e_{i_{1}}^H$ \ is \ 1 \ or \ $e_{i_{1}}$ \ or \ $-1$ \ or \ $-
e_{i_{1}}$.

\bigskip
Second case~: all indices but one are equal, say \ $e_{i_{1}}$.

\medskip
\hskip2cm $v = t_1~e_{i_{1}} + (t_2 + \ldots + t_H)~e_{i_{2}}$
\medskip
\hskip2cm $v^H = \cases{
\big( t_1^2 + (t_2 +\ldots+ t_H)^2\big)^{H/2} &we get \ $0$\cr
\big( t_1^2 + (t_2 +\ldots+ t_H)^2\big)^{(H-1)/2} ~v &we get \
       $e_{i_{1}}/H$\cr
- \big( t_1^2 + (t_2 +\ldots+ t_H)^2\big)^{H/2} &we get \ $0$\cr
- \big( t_1^2 + (t_2 +\ldots+ t_H)^2\big)^{(H-1)/2}~v &we get \ $-
e_{i_{1}}/H$\cr}$
\bigskip
Third case~: at least three different indices
\smallskip
\hskip2cm $v = t_1~e_{i_{1}} + t_2~e_2 + w$

\bigskip
with \ $w$ \ orthogonal to \ $e_{i_{1}}$ \ and \ $e_{i{_2}}$
\medskip
\hskip2cm $v^H = \cases{
( t_1^2 + t_2^2 + w^2)^{H/2}\cr\noalign{\smallskip}
(t_1^2 + t_2^2 + w^2)^{(H-1)/2}~v\cr\noalign{\smallskip}
- (t_1^2 + t_2^2 + w^2)^{H/2}\cr\noalign{\smallskip}
- (t_1^2 + t_2^2 + w^2)^{(H-1)/2}~v\cr}$

\bigskip
We get \ $0$ \ (no homogenous factor in \ $t_1 \ldots t_H$).

\vskip1,5cm
{\bf 2.4  Symmetrization by integral means}
\bigskip\medskip
The main problem of the  $~\hbox{\euro}~$ algorithm is the
disentangling, that is to translate from  $~\hbox{\euro}~$
expression $~\hbox{\euro}~$ \ to an expression without \ $~\hbox{\euro}~ \
~\hbox{\euro}~$ \ using the classical product in the Clifford algebra. 
A tool for that is Dirichlet means, which was studied extensively by
 B.C. Carlson [2] in a completely different situation. He
uses these means for classical special 
functions.

\bigskip
Let $E_{\ell-1}$ be the standard simplex.
\vskip-0,5cm
$$E_{\ell-1} := \{ (t_1,\ldots,t_{\ell-1}) \in \R^{\ell-1}
: ~\forall j, \ t_j\geq 0, \ \sum_{p=1}^{\ell-1} t_p \leq
1\}.$$ 

The beta function in $\ell$ variables is
$$B (b_1,\ldots,b_\ell) := \int_{E_{\ell-1}} ~t_1^{b_1-1}
\ldots t_{\ell-1}^{b_{\ell-1}-1} (1 - t_1 - \ldots -
t_{\ell -1})^{b_\ell-1} \ dt_1 \ldots dt_{\ell-1}$$

$B(b) = B(b_1,\ldots,b_\ell)$  is symmetric. For
$b_j\in\C$, \ $Re ~b_j > 0$  and $g$ integrable, the
Dirichlet measure $\mu_b$ is defined by 
$$\leqalignno{&\int_E g(t) ~d\mu_b (t) :=  &(5)\cr
&\int_{E_{\ell-1}} g
(t_1,\ldots,t_{\ell-1}) {1\over B(b)} \ t_1^{b_1-1}
\ldots t_{\ell-1}^{b_{\ell-1}-1} (1-t_1 - \ldots -
t_{\ell-1})^{b_\ell-1} \ dt_1 \ldots dt_{\ell-1}. &\cr}$$

\bigskip\bigskip
{\it Definition.-}  \ {\it  For $f : S\oplus V
\longrightarrow S\oplus V$ continuous and
$u_1,\ldots,u_\ell$ in $S\oplus V$, put
$$F(f,b,u) := \int_E f(t\!:\!u) ~d\mu_b(t)\leqno (6)$$
with \ $t\!:\! u \ := \displaystyle\sum_{i=1}^{\ell-1}~ t_i u_i
+ \big( 1 - \displaystyle\sum_{i=1}^{\ell-1} t_i\big)
u_\ell$.
}

\bigskip
This integral gives the symmetrization.

A simple illustration with two paravectors $u,v$
$$\eqalign{F (t \rightarrow t^2, 1,1,u,v) &= \int_0^1
(tu+(1-t)v)^2 dt\cr
&= {1\over 3} u^2 + {1\over 3} ~\hbox{\euro}~ uv~\hbox{\euro}~ +
{1\over 3} v^2.\cr}$$

By the remark 1 of paragraph 3, it is always possible to
take only simplices of dimension less than or equal to $n$.

\vskip1,5cm
{\sc 3. -- Analysis with the holomorphic cliffordian product}
\bigskip\bigskip 
{\bf 3.1  Holomorphic cliffordian functions}
\bigskip
In this paragraph, we recall some notions from [9].

Let \ $D$ \ denote the differential operator
$$D = \sum_{i=0}^n \ e_i \ {\partial\over \partial x_i}$$
and let \ $\Delta$ \ be the standard Laplacian
$$\Delta = \sum_{i=0}^n \ {\partial^2\over \partial x_i^2}.$$

If \ $n$ \ is odd, say \ $n = 2m+1$, \ the vector space \ ${\cal V}$ \ of holomorphic
cliffordian functions was defined to be the kernel ot the \ $D\Delta^m$ \ operator.

Let \ $x:= x_0 + \displaystyle\sum_{i=1}^n~e_i x_i$, \ it is holomorphic cliffordian as well
as its powers \ $x^k$ \ (with \ $k\in\Z$). More generally, put \ $\alpha :=
(\alpha_0,\ldots,\alpha_n)$ \ a multiindice, \ $\alpha_i\in\N$, \ and
$$\eqalign{ \vmid{\alpha}
&:= \sum_{i=0}^n~\alpha_i\cr
P_\alpha(x) &:= \sum_{\sigma\in{\hbox{\got S}}} \ \prod_{\nu=1}^{\vmid{~\!\alpha~\!}-1} \
\big( e_{\sigma (\nu)}~ x\big) \ e_{\sigma (\vmid{~\!\alpha~\!})} \cr}$$

where  \ \got S \ is the permutation group with \ $\vmid{\alpha}$ \ elements. By the same
token, put
$$\eqalign{
&\beta := (\beta_0,\ldots,\beta_n), \quad \beta_i\in\N\cr
&\vmid{\beta}~ := \sum_{i=1}^n~Ê\beta_i\cr
&S_\beta (x) := \sum_{\sigma\in\hbox{\got S}}~ \prod_{\nu=1}^{\vmid{~\!\beta~\!}} \big(
x^{-1} \ e_{\sigma (\nu)}\big) x^{-1}.
\cr}$$

The functions \ $P_\alpha$ \ and \ $S_\beta$ \ are, for \ $n$ \ odd, holomorphic cliffordian
but they make sense for all $n$.

Recall from [9] that, when \ $n$ \ is odd there is a Laurent type expansion for holomorphic
cliffordian functions with a pole at the origin~:
$$f(x) = \sum_{\vmid{~\!\beta~\!}~\!Ê<~\! B} S_\beta (x) d_\beta + \sum_{\vmid{~\!\alpha~\!} ~\!\!=
1}^\infty P_\alpha(x) c_\alpha$$
where, in general, \ $d_\beta$ \ and \ $c_\alpha$ \ belong to \ $\R_{0,n}$.

The basic idea is that we work with functions which are limits of sums of \ $x^k$ \ and their
scalar derivatives. Functions generated in this manner are well-defined for all \ $n$. The
problem of building a product is not connected directly with the \ $D\Delta^m$  \ operator.

First we extend the product defined in the previous part.

\bigskip\bigskip\bigskip 
{\bf 3.2  Extension of the product to normally convergent series}
\bigskip
{\sc Theorem 2}.-  Let \ $\displaystyle\sum_{n=0}^\infty a_n$ \ be
a series which converges in norm and such that the coefficients are
products of paravectors. Then the series \
$\displaystyle\sum_{n=0}^\infty ~\hbox{\euro}~ a_n~\hbox{\euro}~$
converges and
$$~\hbox{\euro}~ \sum_{n=0}^\infty a_n~\hbox{\euro}~= \sum_{n=0}^\infty
~\hbox{\euro}~ a_n~\hbox{\euro}~.$$

\bigskip
Proof.- \  From the  inequality (3)
$$\sum_{n=0}^N ~\vvert{~\hbox{\euro}~
a_n~\hbox{\euro}~}~\leq~Ê\sum_{n=0}^N
\vvert{a_n}$$
thus the series \ $\displaystyle\sum_{n=0}^\infty ~\hbox{\euro}~
a_n~\hbox{\euro}~$ \   is convergent in norm.

By linearity
$$~\hbox{\euro}~ \sum_{n=0}^N a_n ~\hbox{\euro}~ - \sum_{n=0}^N
~\hbox{\euro}~ a_n~\hbox{\euro}~ = 0$$
and it suffices to let $N\rightarrow\infty$.

\bigskip
Now it is easy to extend the product to rational functions.
First an example. We define, for $\vvert{1-a}~< 1$
\medskip
\hskip3cm$~\hbox{\euro}~ a^{-1} b~\hbox{\euro}~ :=$
\medskip
\hskip3cm$~\hbox{\euro}~  \big( 1 - (1-a)\big)^{-1} b ~\hbox{\euro}~
=$
\medskip
\hskip3cm$~\hbox{\euro}~\displaystyle\sum_{k=0}^\infty (1-a)^k
b~\hbox{\euro}~ =$

\hskip3cm$\displaystyle\sum_{k=0}^\infty ~\hbox{\euro}~  (1-a)^k
b~\hbox{\euro}~$.

\bigskip
In general we define, for ~$\vvert{1-v_j}~< 1$
\medskip
$$~\hbox{\euro}~ \prod_{i=1}^k~u_i \
\prod_{j=1}^\ell~v_j^{-1}~\hbox{\euro}~ :=
\sum_{k_1=1}^\infty \ldots
\sum_{k_\ell=1}^\infty ~\hbox{\euro}~ \prod_{i=1}^k~u_i~
~\prod_{j=1}^\ell (1-v_j)^{k_j}~\hbox{\euro}~.$$

\bigskip
Of course we have to find the analytic extension for that
symbol. 

A classical example is the following~: for \ $u,v\in (S\oplus V)
\setminus \{ 0\}$
\medskip
$~\hbox{\euro}~ u^{-1} v^{-1}~\hbox{\euro}~$ \ is defined by~: 

if $u$ and $v$
are linearly dependent with \ $v = \lambda u$ \  for some \ 
$\lambda\in\R\setminus\{ 0\}$ \ then it is

$~\hbox{\euro}~ u^{-1} (\lambda u)^{-1} ~\hbox{\euro}~ = \lambda^{-1} ~u^{-2}$.

If \ $u$ \ and \ $v$ \ are linearly independant  for all $t\in [0,1]$, \ $t u+(1-t)v$ \ has an
inverse and we have
$$~\hbox{\euro}~ u^{-1} v^{-1}~\hbox{\euro}~ = \int_0^1 \big(
tu+(1-t)v)^{-2} dt = F(t \fle t^{-1}, 1, 1, u,v).$$
\bigskip
This was introduced in quantum mechanics by R.P. {\sc
Feynmann}~[4].

\bigskip
For a proof, in the open set \ $\vvert{1-u}~< 1$, \
$\vvert{1-v}~< 1$ \   expand in series.

\bigskip
In general, with the hypothesis of linear independence of \ $v_j$
$$~\hbox{\euro}~ \prod_{i=1}^\ell u_i \ \prod_{j=1}^{\ell+1}
v_j^{-1}~\hbox{\euro}~ =$$
$${1\over \ell!} ~\sum_{\sigma\in \hbox{\got S}_\ell}~
\int_E \ \prod_{j=1}^\ell \big( (t\!:\!v)^{-1} u_{\sigma
(j)}\big) \ (t\!:\!v)^{-1} \  dt_1 \ldots dt_\ell.\leqno (7)$$

We have one more \ $v_j$ \ than \ $u_i$. If it is not true, add some \
$v_j = 1$.

\bigskip\bigskip
Remark.- \ Inside the \ $\hbox{\euro}$ \ we compute in the field of
fractions of \ $\R [e_1,\ldots,e_n]$.

\vskip1,5cm
{\bf 3.3 Integral representation formulas for holomorphic
Cliffordian products}
\bigskip\medskip
The standard spectral theory allows us to write
$$f(A) = {1\over 2i\pi} \oint f(z) {1\over z-A}~dz.$$
In particular
$$A^n = {1\over 2i\pi} \oint z^n {1\over z-A}~dz.$$
Now, let $u_1$ and $u_2$ be  linearly
independant elements of the vector space \ $V$,  then
$$~\hbox{\euro}~ u_1^p ~u_2^q ~\hbox{\euro}~ = {1\over (2i\pi)^2}
\oint_{C_1} \oint_{C_2} z_1^p z_2^q \int_0^1 \big( t
(z_1-u_1) + (1-t) (z_2-u_2)\big)^{-2} \ dt~dz_1~dz_2.$$

\medskip
Where \  $C_1$ \  and \  $C_2$ \  are positively oriented
simply closed contours, such that the eigenvalues are inside these
contours.

\bigskip\medskip
For \  $u\in S\oplus V$ \  with \  $u = u_0+\vec u$, the
eigenvalues are \ $u_0 \pm i \ \vvert{\vec u}$

\bigskip
For a general integral representation formula, it is
possible to reduce to the case where $\{
u_1,\ldots,u_\ell\}$ are paravectors and are linearly
independent, then formally~:
\medskip
$(8) \qquad \hbox{\euro}~ f(u_1,\ldots,u_\ell) 
~\hbox{\euro}~ =$
  
$= \displaystyle{1\over (2i\pi)^\ell} \oint_{{\cal C}_1} \ldots
\oint_{{\cal C}_\ell} \ f(z_1,\ldots,z_\ell)~
F(t \rightarrow t^{- \ell},
1,\ldots, 1, z_1-u_1,\ldots, z_\ell-u_\ell)
dz_1 \ldots dz_\ell.$

\vskip1,5cm
{\bf 3.4  Interpolation by polynomials}
\bigskip\medskip
{\sc Theorem 3}.- \ The interpolation formula of Lagrange. Let
$x_0,\ldots,x_\ell$,
\
$\ell +1$ paravectors, \  $a_0,\ldots,a_\ell$, \ $\ell +1$ paravectors.
Put
$$P(x) := \sum_{i=0}^\ell ~\hbox{\euro}~ a_i \prod_{\scriptstyle
k\not= i\atop\scriptstyle k=0}^\ell \ {x - x_k\over x_i -
x_k} \ ~\hbox{\euro}~. \leqno (9)$$
Then, for all $j = 0,\ldots,\ell, \ P(x_j) = a_j$ and, for \ $n$ \ odd, \  $P$
is an holomorphic Cliffordian polynomial of degre $\ell$.

\bigskip\bigskip
Proof.-
$$P(x_j) = \hbox{\euro} ~a_j \prod_{\scriptstyle k\not=
j\atop\scriptstyle k=0}^\ell \ {x_j - x_k\over x_j -
x_k}~Ê\hbox{\euro}~= a_j.$$

The desentangling is easy. Put
$$\leqalignno{
&\alpha_i = \sum_{\scriptstyle k=0\atop\scriptstyle k\not=
i}^\ell ~t_k (x_i-x_k) + t_i + \big( 1 - \sum_{k=0}^\ell
t_k\big)&\cr 
&\beta_{k,i} = \cases{x-x_k &if ~$k\not= i$\cr
                        a_i &if ~$k=i$.\cr}&\cr 
& &\hbox{Then}\cr
&P(x) = \sum_{i=0}^\ell ~ {1\over (\ell +1) !}~
\sum_{\sigma\in\hbox{\got S}_{\ell +1}} ~\int_{E_\ell}
\ \prod_{k=0}^\ell \Big( \alpha_i^{-1} \beta_{\sigma
(k),i}\Big) \alpha_i^{-1}~  dt_0 \  dt_1  ~\ldots ~dt_\ell.&(10)\cr}$$

where \ $\hbox{\got S}_{\ell +1}$ \ is the permutation group of \ $\{
0,1,\ldots,\ell\}$. This formula shows that \ $P$ \ is holomorphic
Cliffordian in \ $x$ \ but also in \ $x_k$ \ and \ $a_k$.

\vskip1,3cm
{\bf 3.5  Product of holomorphic cliffordian functions}
\bigskip\medskip
From the point of view of the product, the \ $S_\beta(x)$ \ are natural~:
\smallskip
put
$$\eqalign{\partial^\beta
&:= {\partial^{\beta_0+\cdots +\beta_n} \over \partial x_0^{\beta_0} \cdots \partial x_n^{\beta_n}
}\cr \noalign{\smallskip}
\hbox{\euro}~ S_\beta (x)~\hbox{\euro}~&= \hbox{\euro}~ (-1)^{\vmid{~\!\beta~\!}}~ \partial^\beta
x^{-1}~\hbox{\euro}~\cr 
&= (-1)^{\vmid{~\!\beta~\!}}~\partial^\beta~ \hbox{\euro}~ x^{-1}~\hbox{\euro}\cr
&= (-1)^{\vmid{~\!\beta~\!}}~\partial^\beta x^{-1}\cr
&= S_\beta (x).\cr}$$

But the \ $P_\alpha(x)$ \ are, in general, different from \ $\hbox{\euro}~
P_\alpha(x)~\hbox{\euro}$. For example~:
$$\hbox{\euro}~ e_1^2x~\hbox{\euro}~ = {1\over 3}~e_1~x~e_1 - {2\over 3}~x.$$
Let
$$k_\alpha := {\vmid{\alpha}~! \over \alpha_0~! ~ \ldots~ \alpha_n~! }$$
we have 
$$\hbox{\euro}~ P_\alpha(x)~\hbox{\euro}~ = k_\alpha \ \partial^\alpha
~x^{2\vmid{~\!\alpha\!~}-1}$$
because the left side is
$$\hbox{\euro}~ P_\alpha(x)~\hbox{\euro}~=~\vmid{\alpha}~! \ ~Ê\hbox{\euro}~
e_0^{\alpha_0} \ldots e_n^{\alpha_n} ~x^{\vmid{~\!\alpha~\!}-1}~\hbox{\euro}$$
and the right side is
$$\eqalign{\partial^\alpha~x^{2\vmid{~\!\alpha~\!}-1}
&= \hbox{\euro}~\partial^\alpha \ x^{2\vmid{~\!\alpha~\!}-1}~\hbox{\euro}\cr 
&=\alpha_0! \ldots \alpha_n! ~\hbox{\euro}~ e_0^{\alpha_0} \ldots e_n^{\alpha_n} \ \
x^{\vmid{~\!\alpha~\!}-1}~\hbox{\euro}.\cr}$$

We may conclude that the set of polynomials \ $\partial^\alpha~x^k$, \ $k\in\N$ \ are better.

For \ $h$ \ and \ $k$ \ in \ $\N$, let
$$\eqalign{
&p(x) = \hbox{\euro}~e_0^{\alpha_0} \cdots e_n^{\alpha_n}~x^h~\hbox{\euro}\cr
&q(x) = \hbox{\euro}~e_0^{\beta_0} \cdots e_n^{\beta_n} \ x^k~\hbox{\euro}.\cr}$$
Then, their product is
$$\hbox{\euro}~p(x) q(x)~\hbox{\euro}~= ~\hbox{\euro}~e_0^{\alpha_0+\beta_0} \cdots
e_n^{\alpha_n+\beta_n} \ x^{h+k}~\hbox{\euro}.$$
\bigskip
Here are other examples of products of holomorphic cliffordian functions.

\bigskip
Product of the exponential and a constant~:
$$\eqalign{\hbox{\euro}~a e^x~\hbox{\euro}~
&= \int_0^1 \ e^{tx} \ a \ e^{(1-t)x}~dt\cr
&= {d\over ds} ~\Big|_{_{\displaystyle s=0}} \ e^{x+sa}.\cr}$$

\bigskip
Product of two exponentials~:
$$\hbox{\euro}~e^x \ e^y~\hbox{\euro}~= \hbox{\euro}~e^{x+y}~\hbox{\euro}~= e^{x+y}.$$

\bigskip
Product of rational functions~:
$$\eqalign{
&\hbox{\euro}~{a\over x-b}~\hbox{\euro} = {{d\over ds}}~Ê\Big|_{_{\displaystyle s=0}} \  \ \int_0^1
\
\big( t + (1-t) (x-b) + sa\big)^{-1} ds\cr 
&\hbox{\euro}~ {1\over (x-a)^p~(x-b)^q}~\hbox{\euro}~ = {(p+q+1)~Ê!\over (p-1)~! \ (q-1)~!} \
\int_0^1 \ \big( ta+(1-t)b\big)^{-(p+q+2)}\  \ t^p (1-t)^q dt.\cr}$$

\bigskip
The computations are the usual ones, by example~:
$$\hbox{\euro}~ {1\over x-a} - {1\over x-b}~\hbox{\euro}~ = ~\hbox{\euro}~ {a-b\over (x-a)
(x-b)}~\hbox{\euro}$$
this means
$$(x-a)^{-1} - (x-b)^{-1} = \int_0^1 \ \big( x-(ta+(1-t)b)\big)^{-1} \ (a-b) \ \big( x -
(ta+(1-t)b)\big)^{-1}~ dt.$$

The basic fact is that  there is no difference between ``variable" and ``constants"~: for \ $n$ \
odd, all expressions are holomorphic cliffordian with respect to their constants too.

\vskip1,3cm
{\bf 3.6  Derivatives and equations of Cauchy-Riemann
type}
\bigskip\medskip
For \ $u\in S\oplus V$, \ $u = \displaystyle\sum_{j=0}^n~u_j e_j$, \
the directional derivative is
$$(u\mid \nabla_x) := \sum_{j=0}^n \ u_j \ {\partial\over\partial
x_j}.$$

\bigskip
Lemma 2.- \  Let $u\in S\oplus V$,
$a\in\R_{0,n} \ 
\ p\in\N$, \ then
$$\leqalignno{(u\mid\nabla_x)  ~\hbox{\euro}~ ax^p ~\hbox{\euro}~ &=
~\hbox{\euro}~ (u\mid\nabla_x) ax^p ~\hbox{\euro}~ &(11)\cr 
&= \cases{0 \  \hbox{if} \  \  p=0\cr
          p  ~\hbox{\euro}~ aux^{p-1} ~\hbox{\euro}~ \ \hbox{if} \  \
p\not= 0.\cr}&\cr}$$

\bigskip
Proof.- \  If \  $p\not= 0$ \  and \  $\varepsilon\in\R$
$$\eqalign{
&(u \mid\nabla_x) ~\hbox{\euro}~ a x^p ~\hbox{\euro}~ = {d\over
d\varepsilon}\Big|_{\varepsilon = 0} ~\hbox{\euro}~ 
a(x+\varepsilon u)^p~\hbox{\euro}~\cr
&= {d\over d\varepsilon}\Big|_{\varepsilon = 0} ~\hbox{\euro}~ a
\sum_{k=0}^p \pmatrix{p\cr k\cr} x^{p-k} \varepsilon^k u^k
~\hbox{\euro}~\cr 
&= {d\over d\varepsilon}\Big|_{\varepsilon = 0}
\ \sum_{k=0}^p \varepsilon^k \pmatrix{p\cr k\cr} ~\hbox{\euro}~ a
x^{p-k} u^k~\hbox{\euro}~\cr 
&= p  ~\hbox{\euro}~ aux^{p-1} ~\hbox{\euro}~.\cr}$$

\bigskip\bigskip
Proposition 2.- \ {\it Let \ $u\in S\oplus V$,
$a\in\R_{0,n} \ \  p\in\Z \setminus \{ 0\}$, \ then
$$(u\mid\nabla_x) ~\hbox{\euro}~ ax^p ~\hbox{\euro}~ = ~\hbox{\euro}~
(u\mid\nabla_x) ax^p ~\hbox{\euro}~ = p ~\hbox{\euro}~ aux^{p-1}
~\hbox{\euro}~.\leqno (12)$$ }

\bigskip
Proof.- \ We have only to work out the case  $p<0$.  If \
$\vvert{1-x}~< 1$
$$\eqalign{
&(u \mid \nabla_x) ~\hbox{\euro}~ a x^{-p} ~\hbox{\euro} \cr
&= (u\mid\nabla_x) ~\hbox{\euro}~a \big( 1 -
(1-x^p)\big)^{-1}~\hbox{\euro}\cr 
&= (u\mid\nabla_x) ~\hbox{\euro}~a
\sum_{q=0}^\infty (1-x^p)^q ~\hbox{\euro}~\cr 
&= \sum_{q=0}^\infty ~\hbox{\euro}~(u\mid\nabla_x) \ a{(1-x^p)}^q
~\hbox{\euro}~\cr 
&= ~\hbox{\euro}~(u\mid\nabla_x) \  a \ \sum_{q=0}^\infty
{(1-x^p)}^q ~\hbox{\euro}~\cr 
&= ~\hbox{\euro}~(u\mid\nabla_x) \  a x^p ~\hbox{\euro}~\cr
&= p ~\hbox{\euro}~ aux^{p-1} ~\hbox{\euro}~.\cr}$$

\bigskip\bigskip
{\sc Theorem 4}.- \ Let $\Omega$ be an open set of $S\oplus V$  with 
$0\in \Omega$. Let  \ 
$f : ~\Omega \rightarrow S\oplus V$   such that locally~:
$$f(x) = \sum_\alpha ~P_\alpha(x) c_\alpha +
\sum_{\vmid{~\beta~}<B} S_\beta (x) d_\beta \leqno (13)$$
with \ $c_\alpha\in\R$ and
$d_\beta\in\R$. Then for all $u\in V$ and $x\not= 0$ \ we have
$${\partial\over\partial x_0} ~\hbox{\euro}~ uf(x) ~\hbox{\euro}~ \ -
(u\mid\nabla_x) \ \hbox{\euro}~ f(x) \ \hbox{\euro} \  = 0. \leqno
(14)$$

\bigskip
Remark.- We get exactly the classical Cauchy-Riemann
equations. When $n = 1$,  that is, in the \ $\C$ \ case, taking $u =
i\lambda$,
$\lambda\in\R$,  we get these well-known equations. When \ $n$ \ is odd, such function is
holomorphic cliffordian and we say that it is with scalar coefficients.

\bigskip\bigskip
Proof.-  By uniform convergence, we have only to compare
$$\eqalign{
&{\partial\over\partial x_0}~\hbox{\euro}~ u P_\alpha (x)
~\hbox{\euro}~= {\partial\over\partial x_0}~\hbox{\euro}~ u k_\alpha
\partial^\alpha \ x^{2\vmid{~\alpha~}-1}~\hbox{\euro}~\cr 
&= k_\alpha ~\partial^\alpha (2\vmid{\alpha}-1)~\hbox{\euro}~
u x^{2\vmid{~\alpha~}-2} ~\hbox{\euro}~\cr  
&(u\mid\nabla_x)~\hbox{\euro}~ P_\alpha (x)~\hbox{\euro}~=~\hbox{\euro}~
(u\mid\nabla_x) k_\alpha \partial^\alpha
x^{2\vmid{~\alpha~}-1} ~\hbox{\euro}~\cr 
&= k_\alpha \partial^\alpha (2\vmid{\alpha}-1)~\hbox{\euro}~ u
~x^{2\vmid{~\alpha~}-2} ~\hbox{\euro}~.\cr}$$

For the \ $S_\beta$, \ we have
$$\eqalign{
&{\partial\over\partial x_0} ~\hbox{\euro}~ u S_\beta (x)
~\hbox{\euro}~ = {\partial\over\partial x_0}~\hbox{\euro}~u h_\beta
~\partial^\beta x^{-1} ~\hbox{\euro}~
= - h_\beta~\partial^\beta ~\hbox{\euro}~ ux^{-2}~\hbox{\euro}~\cr 
&(u\mid\nabla_x) ~\hbox{\euro}~ S_\beta (x)~\hbox{\euro}~=~\hbox{\euro}~
(u\mid\nabla_x) h_\beta~\partial^\beta x^{-1}~\hbox{\euro}~
= h_\beta~\partial^\beta ~\hbox{\euro}~(u\mid \nabla_x) x^{-1}
~\hbox{\euro}~= - h_\beta ~\partial^\beta ~\hbox{\euro}~ u
x^{-2}~\hbox{\euro}~.\cr}$$

\bigskip\bigskip
Remark.- \ For this type of holomorphic Cliffordian function $f$ and for  $x\not= 0$,
$$\lim_{h\rightarrow 0}~\hbox{\euro}~{f(x+h) - f(x)\over h}~
~\hbox{\euro}~,$$
does not depend on the particular paravector $h$, because this is true for $x^p$,  hence also for 
$P_\alpha (x)$, and \ $S_\beta (x)$, and therefore for $f$.

\vskip1,5cm
{\bf 3.7  Taylor formula}
\bigskip\medskip
 Lemma 3.- \ Let \ $p\in\Z$, \ $q\in\N$, \ $u\in V$.  Then
$${\partial^q\over\partial x_0^q} ~\hbox{\euro} ~u^q \ x^p \
\hbox{\euro}~= (u \mid \nabla_x)^q \ x^p.$$

\bigskip
Proof.- iterate (11).

Using scalar derivations this implies
$$\eqalign{
&{\partial^q\over \partial x_0^q} \ \hbox{\euro}~u^q~
P_\alpha(x)~Ê\hbox{\euro}~Ê= (u\mid\nabla_x)^q \ \hbox{\euro}
\ P_\alpha (x) \ \hbox{\euro}\cr 
&{\partial^q\over \partial x_0^q} \
\hbox{\euro}
\ u^q \ S_\beta (x) ~\hbox{\euro}~= (u\mid\nabla_x)^q \
\hbox{\euro}~S_\beta (x)~\hbox{\euro}.\cr}$$

If \ $f$ \ is of the same type as in theorem 4 we have
$${\partial^q\over \partial x_0^q} ~\hbox{\euro}~Êu^q \
f(x)~Ê\hbox{\euro}~= (u\mid \nabla_x)^q \ \hbox{\euro} \  f(x) \
\hbox{\euro}.\leqno (15)$$

\bigskip\bigskip
{\sc Theorem 5 } (Taylor series).- \ Let \ $f$ \ be  an 
holomorphic Cliffordian function with scalar coefficients, then we
have~:
$$\hbox{\euro}~ f(a+x)~\hbox{\euro}~ = \sum_{k=0}^\infty ~{1\over k!} \
\hbox{\euro}~Êx^k
\ {\partial^kf\over
\partial a_0^k} \ (a) ~\hbox{\euro}.$$

\bigskip
Proof.- \ Put \ $x = x_0 + \vec x$.  Since \ $f$ \ is real analytic,
we have
$$\eqalign{ f(a+x)
&= \sum_{k=0}^\infty \ {1\over k!} \ (x\mid\nabla_a)^k \ f(a)\cr
&= \sum_{k=0}^\infty \ {1\over k!} \ \big(
x_0~{\partial\over\partial a_0} + (\vec x\mid\nabla_a)\big)^k \ f(a)\cr
&= \sum_{k=0}^\infty \ {1\over k!} \ \sum_{r+s=k} \ \pmatrix{k\cr
r\cr} \ x_0^r \ {\partial^r\over \partial a_0^r} \ (\vec x \mid
\nabla_a)^s \ f(a)\cr   
\hbox{\euro}~f(a+x) ~Ê\hbox{\euro}~ &= \sum_{k=0}^\infty \ {1\over k!}
\
\sum_{r+s=k}
\
\pmatrix{k\cr r\cr} x_0^r \ {\partial^r\over \partial a_0^r} \
{\partial^s\over
\partial a_0^s} \ \hbox{\euro}~\vec x^s f(a) \ \hbox{\euro}\cr  
&= \sum_{k=0}^\infty \ {1\over k!} \ \sum_{r+s=k} \ \pmatrix{k\cr
r\cr} \ \hbox{\euro} \ x_0^r \ \vec x^s \ {\partial^k\over\partial
a_0^k} \ f(a) \ \hbox{\euro}\cr 
&= \sum_{k=0}^\infty \ {1\over k!} \ \hbox{\euro} \ x^k \
{\partial^k\over \partial a_0^k} \ f(a) \ \hbox{\euro}.\cr}$$

\vskip1cm
{\bf 3.8  Differential calculus}
\bigskip
In this paragraph, \ $n$ \ is odd.

Let $\omega$ be a differential form with values in \ $\R_{0,n}$.  
Then there exist scalar differential forms $\omega_I$ such
that
$$\omega = \Sigma~\omega_I~e_I.$$
We define
$$\hbox{\euro}~\omega~\hbox{\euro}~:=
\Sigma~\omega_I~Ê\hbox{\euro}~e_I~\hbox{\euro}$$
and then the exterior derivative
$$\Eqalign{
&d~\hbox{\euro}~\omega~\hbox{\euro}&&=~
\Sigma~ d\omega_I~\hbox{\euro}~e_I~\hbox{\euro}\cr 
& &&=~\Sigma~\hbox{\euro}~d\omega_I \ e_I~\hbox{\euro}\cr}$$
so that 
$$d~\hbox{\euro}~\omega~\hbox{\euro} =~
\hbox{\euro}~Êd\omega~\hbox{\euro}.$$

\bigskip

Let \ ${\cal P}_v$ \ be the vectorial plane generated by $1$ and
$v$, \ $v\in V$, \ $v^2=-1$.  For a holomorphic Cliffordian
function of the same type as in the previous theorem and
$\Omega_v$ an open set in ${\cal P}_v$ with regular boundary, we have a
Cauchy-Morera theorem.

\bigskip\bigskip
{\sc Theorem 6}.-
$$\eqalign{
&\int_{\partial\Omega_v} \hbox{\euro}~f(x) ~(dx_0+v~d(\vec x
\mid v)~\hbox{\euro}\cr 
&= \int_{\Omega_v} \hbox{\euro}~v {\partial
f(x)\over\partial x_0} - (v\mid\nabla_x)
f(x)~\hbox{\euro}~dx_0 \wedge d(\vec x \mid v)\cr
\noalign{\medskip} 
&= 0.\cr}$$

\bigskip
Proof.-  \ Stokes theorem gives :
$$\eqalign{
&\int_{\partial\Omega_v} \hbox{\euro}~f(x) \big( dx_0+v
d(\vec x \mid v)\big)~\hbox{\euro}\cr 
&= \int_{\Omega_v} d~\hbox{\euro}~f(x) \big( dx_0+v d(\vec x
\mid v)\big)~\hbox{\euro}\cr 
&= \int_{\Omega_v} \hbox{\euro}~df(x) \wedge \big( dx_0 +
v d(\vec x \mid v)\big)~\hbox{\euro}\cr}$$

\bigskip
Then we get
$$\eqalign{
&\int_{\Omega_v} \hbox{\euro} ~(v\mid\nabla) f(x)~d(\vec x
\mid v) \wedge dx_0 + {\partial f(x)\over \partial x_0}~dx_0
\wedge v d(\vec x \mid v)~\hbox{\euro} =\cr
&\int_{\Omega_v} \hbox{\euro} ~v~{\partial f(x)\over\partial
x_0} - (v\mid\nabla) f(x) \hbox{\euro}~dx_0 \wedge d(\vec x
\mid v)  = 0.\cr}$$

\newdimen\margeg \margeg=0pt
\def\bba#1&#2&#3&#4&#5&{\par{\parindent=0pt
    \advance\margeg by 1.1truecm\leftskip=\margeg
    {\everypar{\leftskip=\margeg}\smallbreak\noindent
    \hbox to 0pt{\hss [#1]~~}{{\sc #2} - }#3~; {\it
#4.}\par\medskip
    #5 }
\medskip}}
\vskip1cm
\centerline{References}
\bigskip\bigskip
\bba 1&F.Brackx, R. Delanghe, F. Sommen&Clifford analysis&
        Pitman 1982& &

\bba 2&B.C. Carlson&Special functions of applied Mathematics&
        Academic Press 1977& &

\bba 3&R. Delanghe, F. Sommen, V. Sou\v cek&Clifford
         algebra and spinor-valued functions&Kluwer 1992& &

\bba 4&R.P. Feynman&Space-time approach to quantum
         electrodynamics&Phys. Rev. 76, 769-789, 1949& &

\bba 5&R. Fueter&†ber die analytische Darstellung der
           regulŠren Funktionen einer
           Quaternionenvariablen&Comment. Math. Helv. 8,
           371-378, 1936& &

\bba 6&K. GŸrlebeck, W. Spršssig&Quaternionic and
           Clifford calculus for physicists and engineers&Wiley
           1997& &

\bba 7&D. Hestenes, G. Sobczyk&Clifford algebra to geometric
          calculus&Reidel 1984& &

\bba 8&G. Laville&On Cauchy-Kovalewski extension&Journal of
           functional analysis vol~101, n$^\circ$1, 25-37,
           1991& &

\bba 9&G. Laville, I. Ramadanoff&Holomorphic
           Cliffordian functions&Advances in Clifford algebras
           vol 8, n$^\circ$2, 323-340& &

\bba 10&H. Malonek&Power series representation for monogenic
          functions in \ $\R^{m+1}$ based on a permutational
          product&Complex variables vol 15, 181-191, 1990& &

\bba 11&F. Sommen&A product and an exponential function in
           hypercomplex function theory&Appl. Anal. 12, 13-26 (1981)&&

\bba 12&F. Sommen&The problem of defining abstract bivectors&Result.
           Math. 31, 148-160, (1997)&&

\bba 13&F. Sommen, P. van Lancker&A product for special classes
           of monogenic functions and tensors&Z. Anal. Anwend. 16,
           N$^\circ$4. 1013-1026, (1997)&&

\bba 14&F. Sommen, M. Watkins&Introducing $q$ - Deformation on the
           Level of Vector Variables&Advances in Applied Clifford
           Algebras. Vol 5, n$^\circ$1, 75-82, (1995)&&

\end